\documentclass[12pt]{amsart}

\usepackage{amscd,amssymb,amsmath,amsthm,graphicx,verbatim,fullpage}
\usepackage{textcomp}
\usepackage[OT1,T1]{fontenc}		

\newtheorem{theorem}{Theorem}[section]

\newtheorem{corollary}[theorem]{Corollary}
\newtheorem{proposition}[theorem]{Proposition}

\numberwithin{equation}{section}

\newcommand{\oto}[1]{\stackrel{\underrightarrow{\stackrel{#1}{\,\,\,\,\hphantom{\too}}}}{}} 
\newcommand{\too}{\longrightarrow}
\newcommand{\tens}{\otimes}

\newcommand{\fullwidthfigure}[3]{
	\begin{figure}[htbp]
	\begin{center}\includegraphics[width=10.5cm]{#1}\end{center}
	\caption{#2}
	\label{#3}
	\end{figure}}

\input diagxy

\begin{document}


\title{What separable Frobenius monoidal functors preserve}
\author{Micah Blake McCurdy and Ross Street}
\address{M. B. McCurdy, Macquarie University, NSW 2109, Australia}
\email{micah.mccurdy@gmail.com}
\address{R. Street, Macquarie University, NSW 2109, Australia}
\email{street@math.mq.edu.au}
\label{I1}
\date{April 2009}

\begin{abstract}
\thispagestyle{empty}
\noindent 
Separable Frobenius monoidal functors were defined and studied under that name by Korn\'el Szlach\'anyi~\cite{Szl},~\cite{Szl03}, and by Brian 
Day and Craig Pastro\cite{DP08}. They are a special case of the linearly distributive functors of Robin Cockett and Robert Seely~\cite{CS99}. 
Our purpose is to develop the theory of such functors in a very precise sense. We characterize geometrically which monoidal 
expressions are preserved by these functors (or rather, are stable under conjugation in an obvious sense). We show, by way of corollaries, that
they preserve lax (meaning not necessarily invertible) Yang-Baxter operators, weak Yang-Baxter operators in the 
sense of~\cite{AVR2}, and (in the braided case) weak bimonoids in the sense of~\cite{PS09}. Actually, every weak Yang-Baxter operator is 
the image of a genuine Yang-Baxter operator under a separable Frobenius monoidal functor. Prebimonoidal functors are also 
defined and discussed. 

Les foncteurs mono\"idaux Frobenius s\'eparables ont \'et\'e d\'efinis et \'etudi\'es, sous ce nom, par Korn\'el Szlach\'anyi~\cite{Szl},~\cite{Szl03}, 
et par Brian Day et Craig Pastro~\cite{DP08}. Ils sont un cas sp\'ecial des foncteurs lin\'eaires entre cat\'egories lin\'earement distributives, introduits
par Robin Cockett et Robert Seely~\cite{CS99}. Notre objet est de d\'evelopper la th\'eorie de ces foncteurs en un sens tr\`es pr\'ecis. Nous
caract\'erisons g\'eom\'etriquement les expressions qui sont pr\'eserv\'ees par ces foncteurs (c'est-\`a-dire, sont stables sous conjugaison
en un sens \'evident). Nous montrons sous forme de corollaire qu'ils pr\'eservent les op\'erateurs Yang-Baxter lax (non-n\'ecessairement inversibles),
les op\'erateurs Yang-Baxter faibles dans le sens de~\cite{AVR2}, et (dans le cas tress\'e) les bimono\"ides faibles dans le sens de~\cite{PS09}. 
En fait, chaque op\'erateur Yang-Baxter faible est une image d'un op\'erateur Yang-Baxter v\'eritable par un foncteur Frobenius s\'eparable. Les foncteurs
pr\'ebimono\"daux sont aussi d\'efinis et discut\'es.
\end{abstract}

\maketitle

\vspace*{0.5cm}
Dedicated to Francis Borceux on the occasion of his 60\textsuperscript{th} birthday.


\thispagestyle{empty}
\section{Introduction}

\noindent Frobenius monoidal functors $F:\mathcal{C}\rightarrow \mathcal{X}$
between monoidal categories were defined and studied under that
name in \cite{Szl}, \cite{Szl03} and \cite{DP08} and in a more general
context in \cite{CS99}. If the domain $\mathcal{C}$ is the terminal
category $\text{\boldmath $1$}$, then $F$ amounts to a Frobenius
monoid in $\mathcal{X}$. It was shown in \cite{DP08} that Frobenius
monoidal functors compose, so that, by the last sentence, they take
Frobenius monoids to Frobenius monoids. We concentrate here on separable
Frobenius $F$ and show that various kinds of Yang-Baxter operators
and (in the braided case) weak bimonoids are preserved by $F$.\

We introduce prebimonoidal functors $F:\mathcal{C}\rightarrow \mathcal{X}$
between monoidal categories which are, say, braided. If the domain
$\mathcal{C}$ is the terminal category \nolinebreak $\text{\boldmath $1$}$, then
any (weak) bimonoid in $\mathcal{X}$ gives an example of such an $F$.
We show that prebimonoidal functors compose and relate them to separable
Frobenius functors.\

\section{Definitions}

\noindent Justified by coherence theorems (see \cite{JS} for example), we
write as if our monoidal categories were strict. A functor $F:\mathcal{C}\rightarrow
\mathcal{X}$ between monoidal categories is {\itshape Frobenius}
when it is equipped with a monoidal structure \[ \phi _{A,B}:F A\otimes
F B\rightarrow F( A\otimes B) \qquad \phi _{0}:I\rightarrow F I, \]
and an opmonoidal structure \[ \psi _{A,B}:F( A\otimes B) \rightarrow
F A\otimes F B \qquad \psi _{0}:F I\rightarrow I \] such that

\[ \bfig
\square<1000,500>[F(A \tens B) \tens FC`F(A \tens B \tens C)`FA \tens FB \tens FC`FA \tens F(B \tens C);%
	\phi_{A \tens B,C}`\psi_{A,B} \tens 1`\psi_{A,B \tens C}`1 \tens \phi_{B,C}]
\efig \]\[ \bfig
\square<1000,500>[FA \tens F(B \tens C)`F(A \tens B \tens C)`FA \tens FB \tens FC`F(A \tens B) \tens FC;%
	\phi_{A, B \tens C}`1 \tens \psi_{B,C}`\psi_{A \tens B,C}`\phi_{A,B} \tens 1]
\efig \]

We shall call $F:\mathcal{C}\rightarrow \mathcal{X}$ {\itshape
separable Frobenius monoidal} when it is Frobenius monoidal and
each composite
\[ F( A\otimes B) \overset{\psi _{A,B}}{\longrightarrow }F A\otimes 
F B\overset{\phi _{A,B}}{\longrightarrow }F( A\otimes B) \]
is the identity. We call $F:\mathcal{C}\rightarrow \mathcal{X}$
{\itshape strong monoidal} when it is separable Frobenius monoidal,
$\phi _{A,B}$ is invertible, and $\phi _{0}$ and $\psi _{0}$ are
mutually inverse.\

Suppose $F:\mathcal{C}\rightarrow \mathcal{X}$ is both monoidal
and opmonoidal. By coherence, we have canonical morphisms
\[ \phi _{A_{1}, \ldots , A_{n}}:F A_{1}\otimes \cdots \otimes F
A_{n}\longrightarrow F( A_{1}\otimes \cdots \otimes A_{n}) \]
and
\[
\psi _{A_{1}, \ldots , A_{n}}:F \left( A_{1}\otimes \cdots \otimes
A_{n}\right) \longrightarrow F A_{1}\otimes \cdots \otimes F A_{n}
\]
defined by composites of instances of $\phi $ and $\psi $. If $n=0$
then these reduce to $\phi _{0}$ and $\psi _{0}$; if $n=1$, they
are identities.

The $F${\itshape -conjugate} of a morphism
\[
f: A_{1}\otimes \cdots \otimes A_{n}\longrightarrow B_{1}\otimes
\cdots \otimes B_{m}
\]
in $\mathcal{C}$ is the composite $f^{F}:$
\[ \bfig

\hscalefactor{1.1}
\node a(-1000,+500)[FA_1 \otimes \cdots \otimes FA_n]
\node b(-500,0)[F(A_1\otimes \cdots \otimes A_n)]
\node c(+500,0)[F( B_1\otimes \cdots \otimes B_m)]
\node d(+1000,+500)[FB_1\otimes \cdots \otimes FB_m]
\arrow|a|[a`b;\phi _{A_1,\ldots, A_n}]
\arrow|a|[b`c;Ff]
\arrow|a|[c`d;\psi _{B_1,\ldots,B_m}]

\efig \]

\noindent in $\mathcal{X}$. For $m=1$, this really only requires
$F$ to be monoidal while, for $n=1$, this really only
requires $F$ to be opmonoidal. If a structure in $\mathcal{C}$ is
defined in terms of morphisms between multiple tensors, we can speak
of the $F${\itshape -conjugate} of the structure in $\mathcal{X}$.
For example, we can easily see the well-known fact that the $F$-conjugate
of a monoid, for $F$ monoidal, is a monoid; dually, the $F$-conjugate
of a comonoid, for $F$ opmonoidal, is a comonoid. It was shown in
\cite{DP08} that the $F$-conjugate of a Frobenius monoid is a Frobenius
monoid.

Notice that, for a separable Frobenius monoidal functor $F$, we
have $\phi_{n} \circ \psi_{n} = 1$ for $n>0$.

Suppose $\mathcal{C}$ and $\mathcal{X}$ are braided monoidal. We
say that a separable Frobenius monoidal functor $F: \mathcal C \too \mathcal X$
is {\itshape braided}
when the $F$-conjugate of the braiding $c_{A,B}:A\otimes B\rightarrow
B\otimes A$ in $\mathcal{C}$ is equal to $c_{F A,F B}:F A\otimes
F B\rightarrow F B\otimes F A$ in $\mathcal{X}$. Because of separability,
it follows that $F$ is braided as both a monoidal and opmonoidal
functor.

A {\itshape lax Yang-Baxter (YB) operator }on an object $A$ of a
monoidal category $\mathcal{C}$ is a morphism
$y:A\otimes A\longrightarrow A\otimes A$
satisfying the condition
\[ (y \tens 1) \circ (1 \tens y) \circ (y \tens 1) = (1 \tens y) \circ (y \tens 1) \circ (1 \tens y) \]
A {\itshape Yang-Baxter (YB) operator }is an invertible lax YB-operator.

Recall that the Cauchy (idempotent splitting) completion $\mathcal{Q}\mathcal{C}$
of a category $\mathcal{C}$ is the category whose objects are pairs
$(A,e)$ where $e:A\rightarrow A$ is an idempotent on $A$ and whose
morphisms $f:(A,e)\rightarrow (B,p)$ are morphisms $f:A\rightarrow
B$ in $\mathcal{C}$ satisfying $p f e=f$ (or equivalently $p f=f$
and $f e = f$). Note emphatically that the identity morphism of
$(A,e)$ is $e:(A,e)\rightarrow (A,e)$; in particular, this means
the forgetful $\mathcal{Q}\mathcal{C}\rightarrow \mathcal{C}$, $(A,e)\mapsto
A$, is not a functor. If $\mathcal{C}$ is monoidal then so is
$\mathcal{Q}\mathcal{C}$ with $(A,e)\otimes (A^{\prime },e^{\prime
})=(A\otimes A^{\prime },e\otimes e^{\prime })$ and unit $(I,1)$.\

A {\itshape weak Yang-Baxter operator on} $A$ (compare
\cite{AVR2}) in $\mathcal{C}$ consists of an idempotent $\nabla
:A\otimes A\longrightarrow A\otimes A$, and lax YB-operators $y:A\otimes
A\longrightarrow A\otimes A$ and $y^{\prime }:A\otimes A\longrightarrow
A\otimes A$, subject to the following conditions:
\begin{align}
\nabla \circ y=&y=y\circ \nabla \label{eq21} \\
\nabla \circ y' =& y'=y'\circ \nabla \label{eq22} \\
y\circ y^{\prime }=&\nabla =y^{\prime }\circ y \label{eq23} 
\end{align}
\begin{align}
\left( 1\otimes \nabla \right) \circ \left( \nabla \otimes 1\right) &=
\left( \nabla \otimes 1\right) \circ \left(1\otimes \nabla \right) \\
\left( 1\otimes y\right) \circ \left( \text{$ \nabla $}\otimes 1\right)
&=\left( \nabla \otimes 1\right) \circ \left( 1\otimes y\right), \\
\left( 1\otimes \nabla \right) \circ \left( y\otimes 1\right) &=
\left( y\otimes 1\right) \circ \left( 1\otimes \nabla \right) .
\end{align}
Notice that Equations \nolinebreak \ref{eq21}, \nolinebreak \ref{eq22} \nolinebreak and \nolinebreak \ref{eq23}
say that $y:(A\otimes A,\nabla )\longrightarrow (A\otimes A,\nabla
)$ is a morphism with inverse $y^{\prime }$ in $\mathcal{Q}\mathcal{C}$.

Suppose $(A, \mu:A\otimes A\too A, \eta:I \too A)$ and $(B, \mu :B\otimes B\too A, \eta :I \too B)$ 
are monoids in the monoidal category $\mathcal{C}$. Let a morphism
$\lambda :A\otimes B\longrightarrow B\otimes A$ be given. The following
conditions imply that $A\otimes B$ becomes a monoid with multiplication
$A\otimes B\otimes A\otimes B\overset{1\otimes \lambda \otimes 1}{\longrightarrow
}A\otimes A\otimes B\otimes B\overset{\mu \otimes \mu }{\longrightarrow
}A\otimes B$ and unit $I\overset{\eta \otimes \eta }{\longrightarrow
}A\otimes B$:\
\begin{gather}
\lambda \circ \left( \mu \otimes 1_{B}\right) =\left( 1_{B}\otimes
\mu \right) \circ \left( \lambda \otimes 1_{A}\right) \circ \left(
1_{A}\otimes \lambda \right), \\
\lambda \circ \left( 1_{A}\otimes
\mu \right) =\left( \mu \otimes 1_{A}\right) \circ \left( 1_{B}\otimes
\lambda \right) \circ \left( \lambda \otimes 1_{B}\right) ,
\end{gather}
\begin{equation}
\lambda \circ \left( \eta \otimes 1_{B}\right) =1_{B}\otimes \eta, \qquad
\lambda \circ \left( 1_{A}\otimes \eta \right) =\eta \otimes
1_{A}.%
\label{XRef-Equation-412134328}
\end{equation}
These are the conditions for $\lambda $ to be a {\itshape distributive
law }\cite{Beck69}. A {\itshape weak distributive law} \cite{Stwdl}
is the same except that Equations \nolinebreak \ref{XRef-Equation-412134328}
are replaced by:
\begin{equation}
\left( 1\otimes \mu \right) \circ \left( \lambda \otimes 1\right) \circ
\left( \eta \otimes 1\otimes 1\right) =\left( \mu
\otimes 1\right) \circ \left( 1\otimes \lambda \right) \circ
\left( 1\otimes 1\otimes \eta \right) .
\end{equation}

In the monoidal category $\mathcal{C}$, suppose $A$ is equipped
with a multiplication $\mu :A\otimes A\longrightarrow A$ and a ``switch
morphism'' $\lambda :A\otimes A\longrightarrow A\otimes A$. Supply
$A\otimes A$ with the multiplication $A\otimes A\otimes A\otimes
A\overset{1\otimes \lambda \otimes 1}{\longrightarrow }A\otimes
A\otimes A\otimes A\overset{\mu \otimes \mu }{\longrightarrow }A\otimes
A$. Then a comultiplication $\delta :A\longrightarrow A\otimes A$
preserves multiplication when the following holds:
\begin{equation}
\delta \circ \mu = \left( \mu \otimes \mu \right) \circ \left(
1\otimes \lambda \otimes 1\right) \circ \left( \delta \otimes \delta
\right) .%
\label{multcomult}
\end{equation}
Dually, if we start with $\delta $ and $\lambda $, define the comultiplication
$A\otimes A\overset{\delta \otimes \delta }{\longrightarrow }A\otimes
A\otimes A\otimes A\overset{1\otimes \lambda \otimes 1}{\longrightarrow
}A\otimes A\otimes A\otimes A$ on $A\otimes A$, and ask for $\mu
$ to preserve comultiplication, we are led to the same Equation \nolinebreak
\ref{multcomult}

In a braided monoidal category $\mathcal{C}$, a {\itshape weak bimonoid}
(see \cite{PS09}) is an object $A$ equipped with a monoid structure
and a comonoid structure satisfying Equation \nolinebreak \ref{multcomult}
(with $\lambda =c_{A,A}$) and the ``weak unit and counit'' conditions:
\begin{align}
\varepsilon \circ \mu \circ \left( 1\otimes \mu \right) 
&=\left(\varepsilon \otimes \varepsilon \right) \circ \left( \mu \otimes
\mu \right) \circ \left( 1\otimes \delta \otimes 1\right) \label{weakcounits} \\
&=\left(\varepsilon \otimes \varepsilon \right) \circ \left( \mu \otimes
\mu \right) \circ \left( 1\otimes c_{A,A}^{-1}\otimes 1\right) \circ
\left( 1\otimes \delta \otimes 1\right) \nonumber \\
\left( 1\otimes \delta \right) \circ \delta \circ \eta 
&=\left( 1\otimes \mu \otimes 1\right) \circ \left( \delta \otimes \delta \right)
\circ \left( \eta \otimes \eta \right) \label{weakunits} \\
& =\left( 1\otimes \mu \otimes 1\right) \circ \left( 1\otimes c_{A,A}^{-1}\otimes 1\right) \circ
\left( \delta \otimes \delta \right) \circ \left( \eta \otimes \eta \right) \nonumber
\end{align}

A {\itshape lax Yang-Baxter (YB) operator }on a functor $T: \mathcal{A}\rightarrow
\mathcal{C}$ into a monoidal category $\mathcal{C}$ is a natural
family of morphisms
\[ y_{A,B}:T A\otimes T B\longrightarrow T B\otimes T A \]
satisfying the condition
\[ \bfig

\hscalefactor{1.1}
\vscalefactor{0.7}
\node l(-1000,0)[TA \tens TB \tens TC]
\node tl(-500,+500)[TB \tens TA \tens TC]
\node tr(+500,+500)[TB \tens TC \tens TA]
\node r(+1000,0)[TC \tens TB \tens TA]
\node bl(-500,-500)[TA \tens TC \tens TB]
\node br(+500,-500)[TC \tens TA \tens TB]

\arrow|a|[l`tl;y \tens 1]
\arrow|a|[tl`tr;1 \tens y]
\arrow|a|[tr`r;y \tens 1]
\arrow|b|[l`bl;1 \tens y]
\arrow|b|[bl`br;y \tens 1]
\arrow|b|[br`r;1 \tens y]

\efig \]

One special case is where $\mathcal{A}=1$ so that $T$ is an object
of $\mathcal{C}$: then we obtain a lax YB-operator on the object
$T$ as above. Another case is where $\mathcal{A}=\mathcal{C}$ and
$T$ is the identity functor: each (lax) braiding $c$ on $\mathcal{C}$
gives an example with $y_{A,B}=c_{A,B}$.

Suppose $T: \mathcal{A}\rightarrow \mathcal{C}$ is a functor and
$F:\mathcal{C}\rightarrow \mathcal{X}$ is a functor between monoidal
categories. Suppose lax YB-operators $y$ on $T$ and $z$ on $F T$
are given. We define $F$ to be {\itshape prebimonoidal} {\itshape
relative to} $y$ {\itshape and} $z$ when it is monoidal and opmonoidal,
and satisfies

\[ \bfig

\hscalefactor{1.09}
\node l(-800,0)[F(TA \tens TB) \tens F(TC \tens TD)]
\node tl(-650,+500)[FTA \tens FTB \tens FTC \tens FTD]
\node tr(+650,+500)[FTA \tens FTC \tens FTB \tens FTD]
\node r(+800,0)[F(TA \tens TC) \tens F(TB \tens TD)]
\node bl(-650,-500)[F(TA \tens TB \tens TC \tens TD)]
\node br(+650,-500)[F(TA \tens TC \tens TB \tens TD)]

\arrow|a|[l`tl;\psi \tens \psi]
\arrow|a|[tl`tr;1 \tens z \tens 1]
\arrow|a|[tr`r;\phi \tens \phi]
\arrow|b|[bl`br;F(1 \tens y \tens 1)]
\arrow|b|[l`bl;\phi]
\arrow|b|[br`r;\psi]

\efig \]

When $\mathcal{C}$ and $\mathcal{X}$ are (lax) braided and $T$ is
the identity with $y_{A,B}=c_{A,B}$ and $z_{A,B}=c_{F A,F B}$, we
merely say $F$ is {\itshape prebimonoidal}. Such an $F$ is bimonoidal when,
furthermore, $FI$, with its natural monoid and comonoid structure, is a bimonoid.
We were surprised not to find this concept in the literature, however, we have found
that it was presented in preliminary versions of the forthcoming book \cite{AM}, and
in talks by the authors of the same.

\section{Separable invariance and connectivity}

\begin{figure}[htbp]
\begin{center}
\resizebox{300pt}{!}{\includegraphics{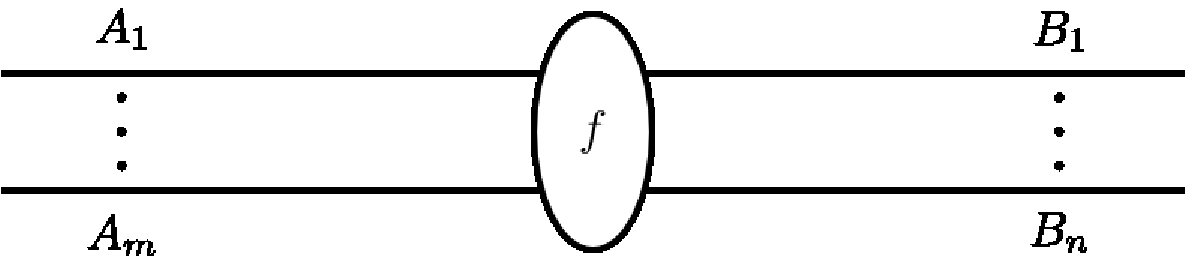}}
\end{center}
\caption{}
\label{XRef-Figure-1013172236}
\end{figure}

We begin by reviewing some concepts from \cite{JS2}. Progressive
plane string diagrams are deformation classes of progressive plane
graphs. Here we will draw them progressing from left to right (direction
of the $x$-axis) rather than from down to up (direction of the $y$-axis).
A tensor scheme is a combinatorial directed graph with vertices
and edges such that the source and target of each edge is a word
of vertices (rather than a single vertex). Progressive string
diagrams $\Gamma $ can be labelled (or can have valuations) in a
tensor scheme $\mathcal{D}$: for a given labelling $v:\Gamma \rightarrow
\mathcal{D}$, the labels on the edges (strings) $\gamma $ of $\Gamma
$ are vertices $v( \gamma ) $ of $\mathcal{D}$ while the labels
on the vertices (nodes) $x$ of $\Gamma $ are edges $v( x) :v( \gamma
_{1}) \cdots v( \gamma _{m}) \rightarrow v( \delta _{1}) \cdots
v( \delta _{n}) $ of $\mathcal{D}$ where $\gamma _{1}, \cdots ,\gamma
_{m}$ are the input edges and $\delta_1, \cdots, \delta_n$
are the output edges of $x$ read from top to bottom; see Figure \nolinebreak
\ref{XRef-Figure-1013172236} where $f=v( x) $, $A_{1}=v( \gamma
_{1}) $, $B_{n}=v( \delta _{n}) $, and so on. The free monoidal
category $\mathcal{F}\mathcal{D}$ on a tensor scheme $\mathcal{D}$
has objects words of vertices and morphisms progressive plane string
diagrams labelled in $\mathcal{D}$; composition progresses horizontally
while tensoring is defined by stacking diagrams vertically.

\begin{figure}[htbp]
\begin{center}
\resizebox{200pt}{!}{\includegraphics{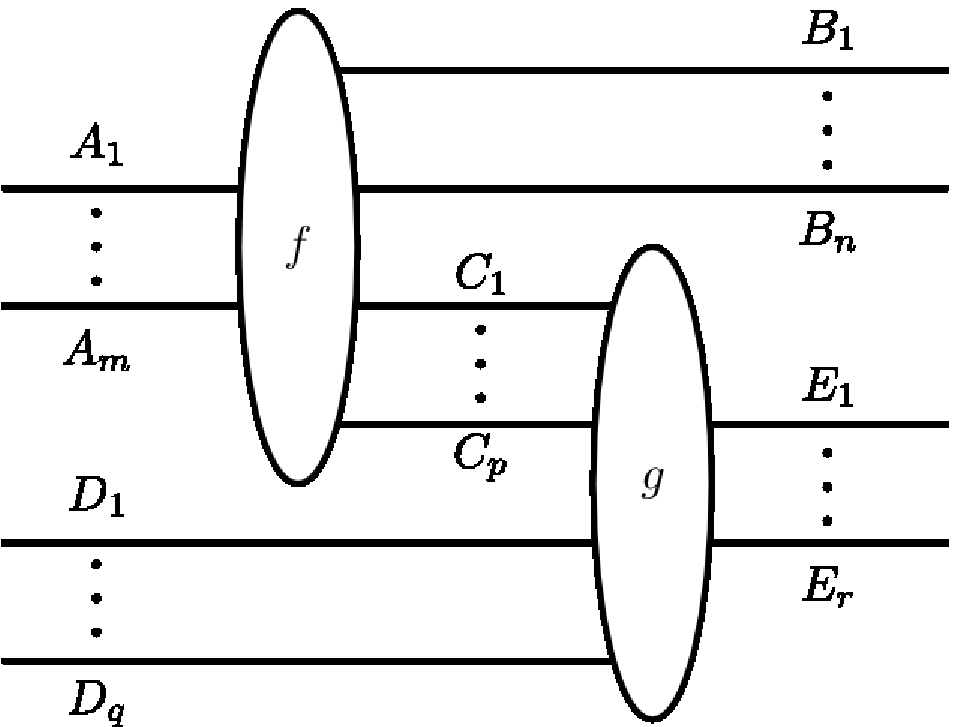}}
\end{center}
\caption{}
\label{cut}
\end{figure}

Every monoidal category $\mathcal{C}$ has an underlying tensor scheme:
the vertices are the objects of $\mathcal{C}$ and the edges from
one word $A_1 \cdots A_{m}$ of objects to another $B_{1} \cdots
B_{n}$ is a morphism $f:A_{1} \otimes \cdots \otimes A_{m}\rightarrow
B_{1}\otimes \cdots \otimes B_{n}$ in $\mathcal{C}$; see Figure \nolinebreak
\ref{XRef-Figure-1013172236}. When we speak of a labelling of a
string diagram in $\mathcal{C}$ we mean a labelling in the underlying
tensor scheme; here we will simply call this a {\itshape string
diagram in} $\mathcal{C}$. The value $v( \Gamma ) $ of the string
diagram $v:\Gamma \rightarrow \mathcal{C}$ is a morphism obtained
by deforming $\Gamma $ so that no two vertices of $\Gamma $ are
on the same vertical line then by horizontally composing strips
of the form
\[
1_{C_{1}}\otimes \cdots \otimes 1_{C_{h}}\otimes f\otimes 1_{D_{1}}\otimes
\cdots \otimes 1_{D_{k}}.
\]
Calculations in monoidal categories can be performed using string
diagrams rather than the traditional diagrams of category theory.
The value of Figure \nolinebreak \ref{XRef-Figure-1013172236} is of course $f$.
Figure \nolinebreak \ref{cut} shows a string diagram $v:\Gamma
\rightarrow \mathcal{C}$ whose value $v( \Gamma ) $ is
\[ \bfig

\node a(0,+500)[A_{1} \otimes \cdots \otimes A_{m}\otimes D_{1} \otimes \cdots \otimes D_{q}]
\node b(0,0)[B_1 \otimes \cdots \otimes B_n \otimes C_{1} \otimes \cdots \otimes C_{p}\otimes D_{1}\otimes \cdots \otimes D_{q}]
\node c(0,-500)[E_{1} \otimes \cdots \otimes E_{p}]
\arrow|m|[a`b;f \otimes 1]
\arrow|m|[b`c;1 \otimes g]

\efig \]

Now we return to our study of separable Frobenius monoidal functors.

\begin{figure}[htbp]
\begin{center}
\resizebox{200pt}{!}{\includegraphics{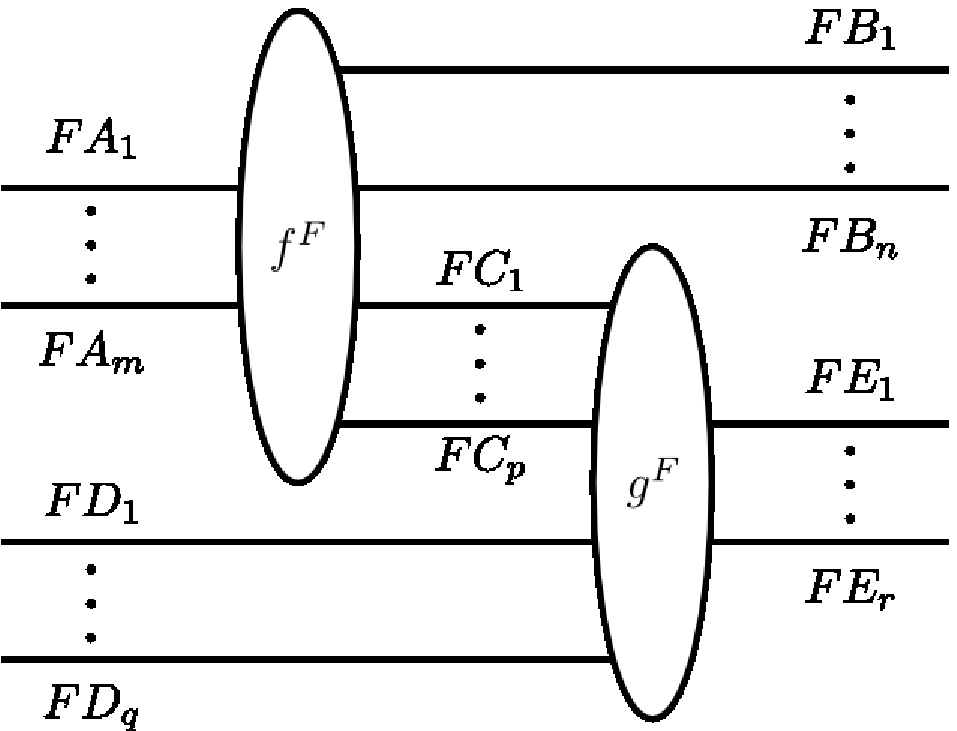}}
\end{center}
\caption{}
\label{XRef-Figure-101319744}
\end{figure}

Suppose $v:\Gamma \rightarrow \mathcal{C}$ is a string diagram in
a monoidal category $\mathcal{C}$ and $F:\mathcal{C}\rightarrow
\mathcal{X}$ is a monoidal and opmonoidal functor. We obtain a {\itshape
conjugate string diagram} $v^{F}:\Gamma \rightarrow \mathcal{X}$
in $\mathcal{X}$ by defining\
\[ v^F(\gamma) = Fv( \gamma ) \mbox{ and } v^F(x) = {v(x)}^F \]
for each edge $\gamma $ and each node $x$ of $\Gamma $. The conjugate
of the string diagram in Figure \nolinebreak \ref{cut} is
shown in Figure \nolinebreak \ref{XRef-Figure-101319744}.
A (progressive plane) string diagram $\Gamma $ is called {\itshape
[separable] Frobenius invariant} when, for any labelling $v:\Gamma
\rightarrow \mathcal{C}$ of $\Gamma $ in any monoidal category $\mathcal{C}$
and any [separable] Frobenius monoidal functor $F:\mathcal{C}\rightarrow
\mathcal{X}$, the value of the conjugate diagram $v^{F}$ in $\mathcal{X}$
is equal to the conjugate of the value of $v$; that is,
\begin{equation}
v^{F}( \Gamma ) ={v( \Gamma ) }^{F}.%
\label{XRef-Equation-1013195041}
\end{equation}


As mentioned before, for a separable Frobenius monoidal functor
$F$, we have $\phi_n \circ \psi_n = 1$ for $n>0$.

The following two theorems characterize which string diagrams are preserved by Frobenius and separable Frobenius
monoidal functors in terms of connectedness and acyclicity. Robin Cockett pointed out to us that similar geometric 
conditions occur in the work of Girard \cite{Girard} and Fleury and Retor\'e (\cite{FleuryRetore}, \S 3.1).
There may be a relationship with our results but the precise nature is unclear.

\begin{theorem}
A progressive plane string diagram is separable Frobenius invariant if and only if it is connected.
\label{septheorem}
\end{theorem}
\begin{proof}

In Figure \nolinebreak \ref{longstrings}, we show that Equation \nolinebreak \ref{XRef-Equation-1013195041}
holds for the string diagram $v:\Gamma \rightarrow \mathcal{C}$
as in Figure \nolinebreak \ref{cut}, provided $p>0$ (as required
for $\Gamma $ to be connected). To simplify notation we write
$F A$ for $F( A_{1} \otimes \cdots \otimes A_{m}) $ and write $A^{F}$
for $F A_{1} \otimes \cdots \otimes F A_{m}$. We also leave out
some tensor symbols $\otimes $. The second equality in Figure \nolinebreak \ref{longstrings}
is where separability, and the fact that the length $p$ of the word
$C$ is strictly positive, are used; the third is where a Frobenius
property is used.

\fullwidthfigure{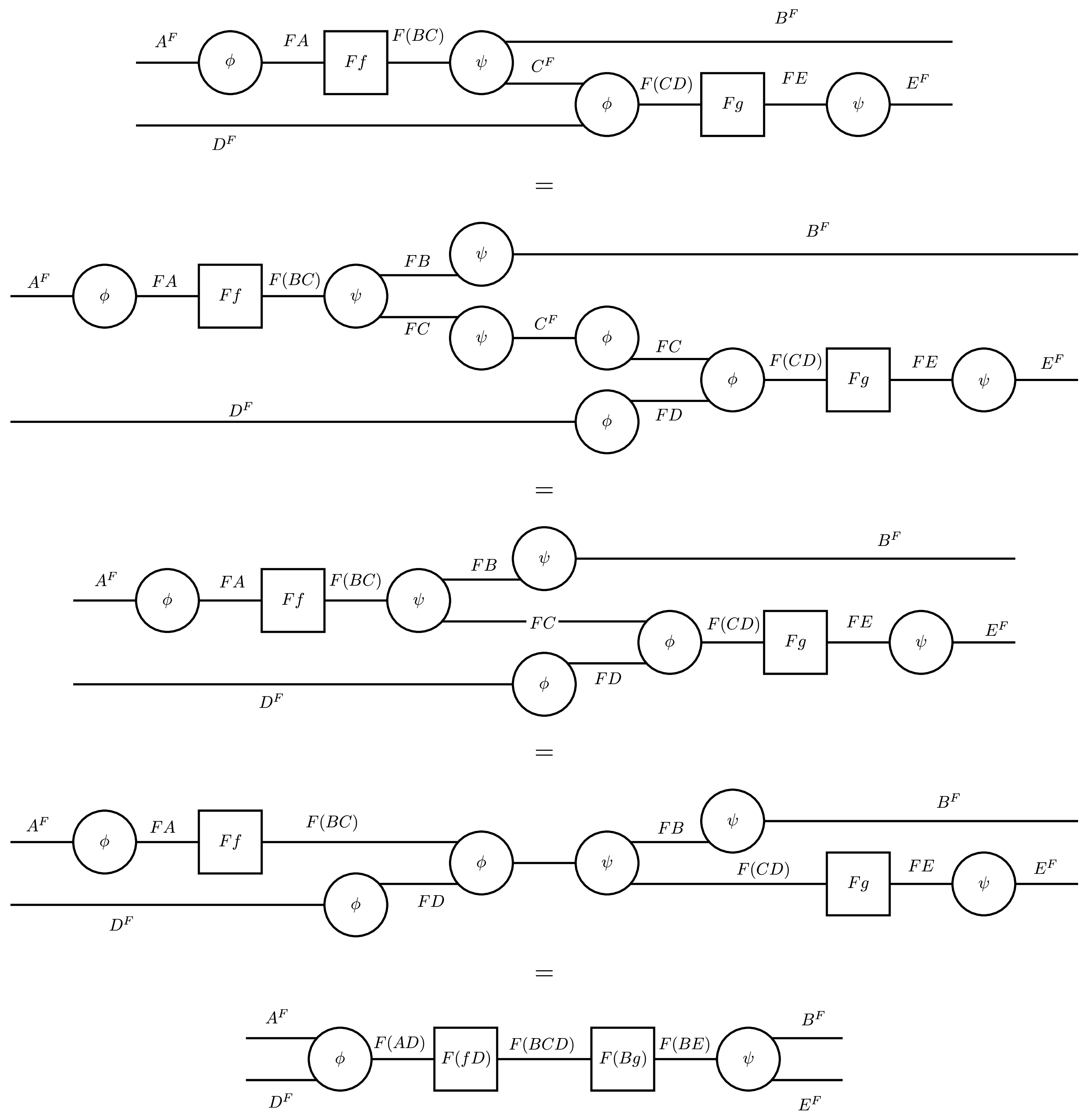}{}{longstrings}


Similarly, an obvious (horizontal) dual diagram to Figure \nolinebreak \ref{cut}
(look through the back of the page!) can be shown separably invariant.
Furthermore, it is simple to show that diagrams of the form shown in Figure \nolinebreak \ref{internalcut}
are separably invariant, as well as their diagrams of the dual form.
\begin{figure}[htbp]
\begin{center}
\resizebox{200pt}{!}{\includegraphics{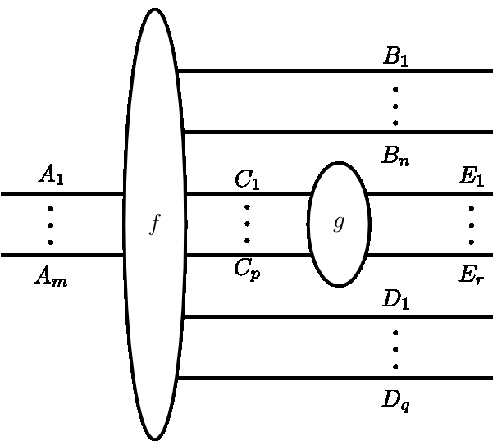}}
\end{center}
\caption{}
\label{internalcut}
\end{figure}

By a similar proof to the above, such diagrams and their duals are separably invariant.
Every connected string diagram can be constructed by iterating these
four processes, this proves ``if''. For ``only if'', we exploit the fact that every
string diagram can be interpreted in the terminal monoidal category $1$ and that separable
Frobenius monoidal functors $F: 1 \rightarrow \mathcal C$ are precisely separable Frobenius
algebras in $\mathcal C$. 

Suppose for a contradiction that a disconnected string diagram $\Gamma$ with
$n$ input wires and $m$ output wires is invariant under conjugation by such a separable
Frobenius $F$, that is, a separable Frobenius algebra $C$. This asserts
the equality of two morphisms $C^{\tens n} \too C^{\tens m}$, the first (obtained by
taking the (trivial) value of the labelling in $1$ and then applying $F$) is the composite
of $n$-fold multiplication followed by $m$-fold comultiplication; the second (obtained by
applying $F$ to the labelling and then taking the value in $\mathcal C$) is considerably
more complicated, containing at least two connected components since $\Gamma$ is assumed to
be disconnected.
By prepending $n$ units and appending $m$ counits, the first becomes the barbell
of unit followed by counit; the latter becomes an endomap of the tensor unit of $\mathcal C$
with at least two connected components. If the tensor product of $\mathcal C$ is symmetric,
this last simplifies to as many copies of the barbell as there are connected
components of $\Gamma$; hence, it suffices to find a separable Frobenius algebra for which
the barbell does not equal any non-trivial power of itself.

We give two examples of such separable Frobenius algebras, a simple
algebraic example and a more complicated geometric example. First, consider the complex numbers
as a Frobenius algebra over the reals. Kock (\cite{Kock}, Example 2.2.14) notes that $\mathbb C \too \mathbb R$ given
by $x+iy \mapsto ax + by$ is a Frobenius form on $\mathbb C$ for all $a$ and $b$ not both zero. Choosing $a = 2, b = 0$
gives a \emph{separable} Frobenius structure, and the ``barbell'' 
$\mathbb R \too \mathbb C \too \mathbb R$ is multiplication by 2, which does not equal any non-trivial
power of itself. This completes the proof of the converse of the theorem.

We sketch the construction of a more complicated but perhaps more pleasing, geometric example: 
consider the category
$\bf 2Thick$, as described in \cite{Lauda}, 
whose objects are finite disjoint unions of the interval (identified with the natural numbers), embedded
in the plane, and whose morphisms are 
boundary-preserving-diffeomorphism classes smooth oriented surfaces embedded
in the plane with boundary equal to the union of domain and codomain. For instance, 
Figure \nolinebreak \ref{pants} shows a morphism in $\bf 2Thick$ from 2 to 1.
\begin{figure}[htbp]\begin{center}\includegraphics[width=5cm]{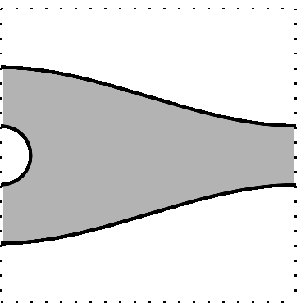}\end{center}\caption{A morphism in $\bf 2Thick$}\label{pants}\end{figure}
Lauda proves that $\bf 2Thick$ is the free monoidal category containing a Frobenius algebra; the morphism from 2 to 1 shown in 
Figure \nolinebreak \ref{pants}
is the multiplication for this Frobenius algebra; the obvious similar map from 1 to 2 is the comultiplication.
However, for this theorem, we require a \emph{separable} Frobenius algebra, so we modify $\bf 2Thick$ to obtain a category in which the equality
in Figure \nolinebreak \ref{bubble} holds; in fact, we conjecture, to obtain the free monoidal category containing a separable Frobenius algebra.
\fullwidthfigure{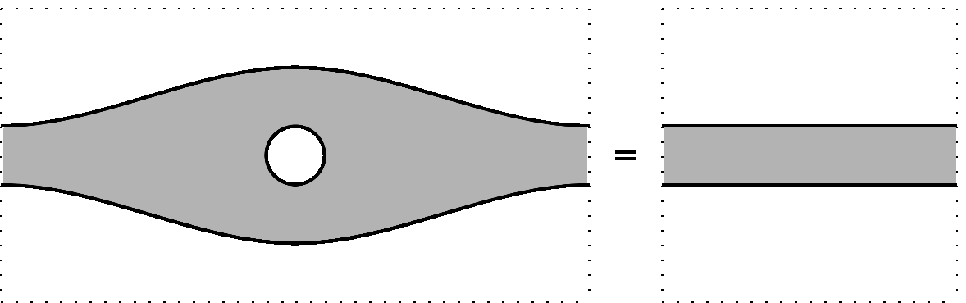}{}{bubble}
Specifically, instead of taking boundary-preserving diffeomorphism classes of morphisms, we say that two morphisms $k \too l$ are equal if there is a
suitable 3-manifold $M$ with corners which can be embedded in the unit cube in such that the intersection with the top face is the first morphism and
the intersection with the bottom face is the second morphism. Here ``suitable'' means that the 3-manifold must be trivial on the domains and codomains
$k$ and $l$, and, crucially, the only critical points of the boundary of $M$ permitted are ``cups'' -- that is, critical points which are not
saddle points where the convex portion of the critical point lies \emph{outside} the manifold $M$. There is an evident such ``cup'' which will witness the 
desired equality shown in Figure \nolinebreak \ref{bubble}. Let us call this quotient of $\bf 2Thick$ by the name $\bf 2Thick'$.

Most importantly, it is clear that no two morphisms with different numbers
of connected components can be identified by this equivalence relation, so any disconnected string diagram will fail to be separably invariant
with respect to the canonical separable Frobenius functor $1 \too {\bf 2Thick'}$. 						\end{proof}

What this implies is that separable Frobenius monoidal functors
preserve equations in monoidal categories for which both sides of
the equation are values of connected string diagrams. For example:
\begin{corollary}

For $n>1$, equations of the form:
\[ \left( a_n\otimes 1\right) \left( 1\otimes a_{n-1}\right)
\left( a_{n-2}\otimes 1\right) \cdots = \left( 1\otimes
b_n\right) \left( b_{n-1}\otimes 1\right) \left( 1\otimes
b_{n-2}\right) \cdots 
, \] involving morphisms
\[ a_1, \ldots, a_n, b_1, \ldots, b_n:A\otimes A\longrightarrow A\otimes A ,\]
are stable under F-conjugation. In fact, for $n=2$, Frobenius $F$ will do.
\end{corollary}

The proof of Theorem \nolinebreak \ref{septheorem} can be slightly modified to give the analagous 
result for merely Frobenius monoidal functors instead of separable Frobenius monoidal
functors.
\begin{theorem}
A progressive plane string diagram is Frobenius invariant if and only if it is connected and simply connected.
\end{theorem}
\rm 
\begin{proof} We have noted that all connected string diagrams can be obtained as iterations of the constructions shown in
Figures \nolinebreak \ref{cut} \nolinebreak and \nolinebreak \ref{internalcut} and their duals, with the restriction that $p > 0$. 
All connected and simply connected string diagrams can be obtained in this way with
the restriction that $p=1$. The only step of the proof in Figure \nolinebreak \ref{longstrings} (and also the corresponding proof for
the case shown in Figure \nolinebreak \ref{internalcut}) which requires separability is
the cancellation of $FC \oto{\psi} C^F \oto{\phi} FC$ to obtain the identity on $FC$; since $p = 1$, we have $FC = C^F = FC_1$, and both
of these maps are identities. Hence, the same proof will go through in this case, establishing ``if''. 

Conversely, suppose that $\Gamma$ is a string diagram which is not connected and acyclic. By Theorem \nolinebreak \ref{septheorem}, we may assume
that $\Gamma$ is connected and therefore is not acyclic. Then for $\Gamma$ to be invariant under the canonical Frobenius monoidal functor
$1 \too {\bf 2Thick}$ described in \cite{Lauda} and referred to already in Theorem \nolinebreak \ref{septheorem} would imply that there is a
diffeomorphism between two 2-manifolds of different genus; this is not the case.		\end{proof}

\begin{corollary}
Weak bimonoids are preserved by braided separable Frobenius functors.
\end{corollary}
\begin{proof}
Weak bimonoids satisfy Equations \nolinebreak \ref{multcomult}, \nolinebreak \ref{weakcounits}, \nolinebreak and \nolinebreak \ref{weakunits} 
these equations are labelled versions of the following string-diagrams:
\begin{center}\includegraphics[width=\textwidth]{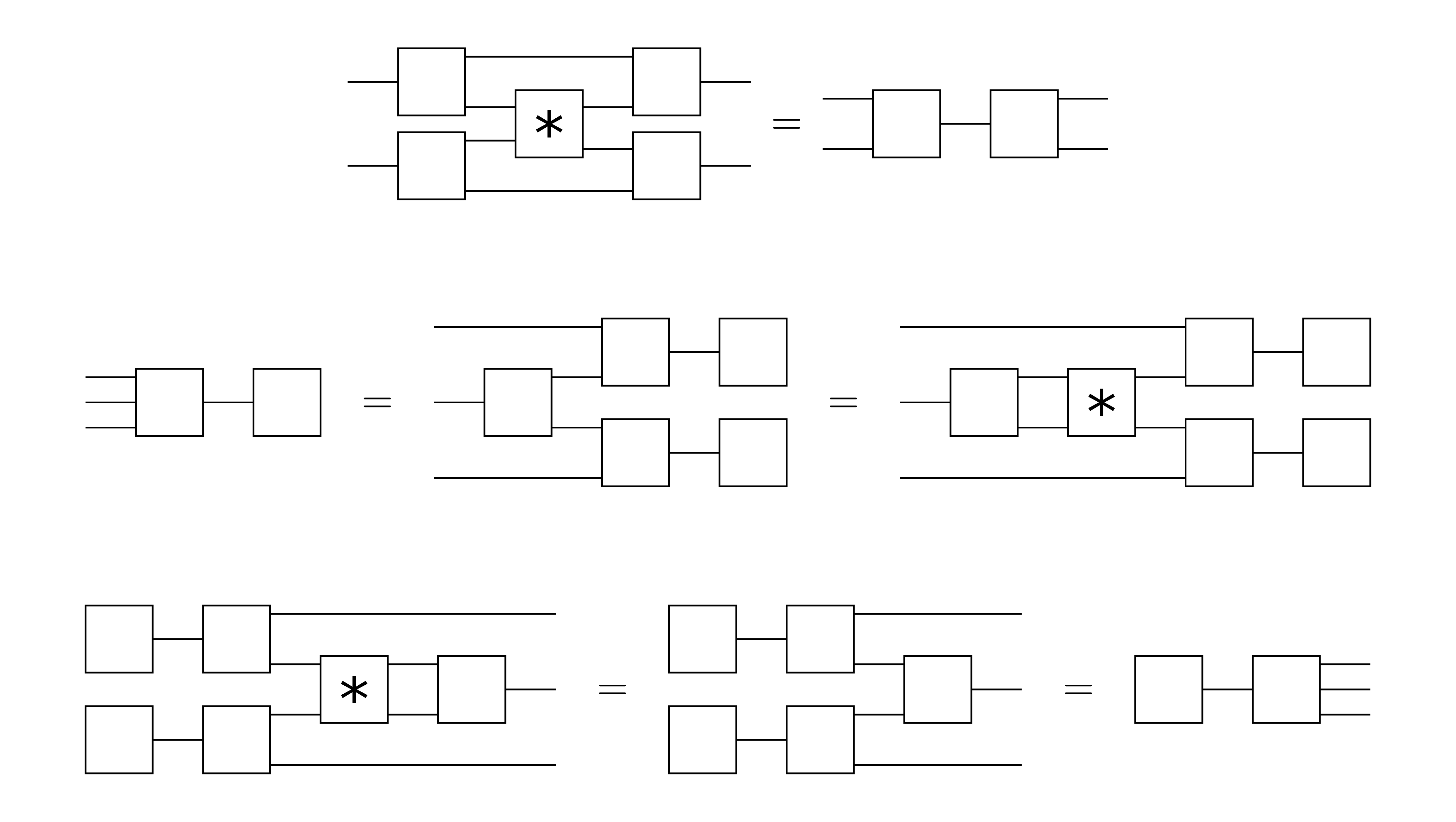}\end{center}
These are clearly connected and hence preserved by separable Frobenius functors.
The asterisks indicate labellings by braids or their inverses, these are preserved by braided Frobenius functors.
\end{proof}

However, genuine bimonoids are not preserved in general: the three unit
and counit equations for a bimonoid involve non-connected
string diagrams.

\begin{corollary}

Weak distributive laws are preserved under $F$-conjugation.
\end{corollary}

However, distributive laws are not preserved in general: the string
diagrams for the right-hand sides of Equations \nolinebreak \ref{XRef-Equation-412134328}
are not connected.

\begin{corollary}
Lax YB-operators are preserved under $F$-conjugation.
\end{corollary}

However, YB-operators are not preserved: invertibility involves
an equation whose underlying diagram is a pair of disjoint strings
and so is disconnected.

\begin{corollary}
Weak YB-operators are preserved under $F$-conjugation. In particular,
the $F$-conjugate of a YB-operator is a weak YB-operator.
\end{corollary}

\begin{proposition}
Every weak YB-operator in a monoidal category in which idempotents split is the conjugate of a YB-operator under some
separable Frobenius monoidal functor.
\end{proposition}
{\itshape Proof. }Let $\mathcal{C}$ be a such a monoidal category containing
an object $D$ and an idempotent $\nabla : D\otimes D\longrightarrow
D\otimes D$ such that $(\nabla \otimes 1)(1\otimes \nabla )=(1\otimes
\nabla )(\nabla \otimes 1)$. Then there is an idempotent $\nabla_n : D^{\otimes n}\longrightarrow D^{\otimes n}$ recursively
defined by: \begin{align*}
	\nabla_0 &= 1_I \\
	\nabla_1 &= 1_D \\
	\nabla_2 &= \nabla \\
	\nabla_n &= (1 \tens \nabla_{n-1}) \circ (\nabla \tens 1) \qquad \mbox{for $n >2$}
\end{align*}
\noindent Let $\mathcal{C}( D) $ be the subcategory of $\mathcal{Q}\mathcal{C}$
whose objects are the pairs $(D^{\otimes n},\nabla _{n})$ and whose
morphisms $f:(D^{\otimes n},\nabla _{n})\rightarrow (D^{\otimes
m},\nabla _{m})$ are those in $\mathcal{Q}\mathcal{C}$ for which:
\begin{gather}
(1\otimes f)(\nabla _{n}\otimes 1)=(\nabla _{m}\otimes 1)(1\otimes f) \\
(f\otimes 1)(1\otimes \nabla _{n})=(1\otimes \nabla _{m})(f\otimes 1).
\end{gather}
The category $\mathcal{C}( D) $ becomes monoidal via 
\[ (D^{\otimes n},\nabla_n)\otimes (D^{\otimes m},\nabla _{m})=(D^{\otimes(n+m)},\nabla_{n+m}). \]
Note that this is not the same as the usual tensor product on $\mathcal{QC}$ which is inherited
from that of $\mathcal C$. A weak Yang-Baxter operator on $D$ in
$\mathcal{C}$ is a Yang-Baxter operator on $D$ in $\mathcal{C}(
D) $. Since idempotents split in $\mathcal{C}$ then we have a functor $\mathcal{C}(
D) \rightarrow \mathcal{C}$ taking each idempotent to a splitting.
Moreover, this functor $\mathcal{C}( D) \rightarrow \mathcal{C}$ is separable
Frobenius (although not strong) and so each weak YB-operator is the image of a genuine YB-operator.

\begin{proposition} \label{pfc}
Prebimonoidal functors compose.
\end{proposition}
\begin{proof}
Suppose that $F: \mathcal C \too \mathcal X$ is prebimonoidal with respect to a YB-operator $y$ on $T : \mathcal A \too \mathcal C$ and 
a YB-operator $z$ on $FT$, and suppose further that $G : \mathcal X \too \mathcal Y$ is prebimonoidal with respect to $z$ and a YB-operator
$a$ on $GFT$. Then the diagram in Figure \nolinebreak \ref{lozenge} proves that $GF$ is prebimonoidal with respect to $y$ and $a$.

\begin{figure}[htbp]
\[ \bfig

\hscalefactor{0.8}
\vscalefactor{0.5}
\node a(0,+3000)[GF(Tu \tens Tv) \tens GF(Tw \tens Tx)]
\node b(+1000,+2000)[G(FTu \tens FTv) \tens G(FTw \tens FTx)]
\node c(+1000,+1000)[GFTu \tens GFTv \tens GFTw \tens GFTx]
\node d(+1000,-1000)[GFTu \tens GFTw \tens GFTv \tens GFTx]
\node e(+1000,-2000)[G(FTu \tens FTw) \tens G(FTv \tens FTx)]
\node f(0,-3000)[GFTu \tens GFTw \tens GFTv \tens GFTx]
\node g(-1000,+2000)[G(F(Tu \tens Tv) \tens F(Tw \tens Tx))]
\node h(-1000,+1000)[GF(Tu \tens Tv \tens Tw \tens Tx)]
\node i(-1000,-1000)[GF(Tu \tens Tw \tens Tv \tens Tx)]
\node j(-1000,-2000)[G(F(Tu \tens Tw) \tens F(Tv \tens Tx))]
\node k(0,+500)[G(FTu \tens FTv \tens FTw \tens FTx)]
\node l(0,-500)[G(FTu \tens FTw \tens FTv \tens FTx)]

\arrow|r|[a`b;G\psi \tens G\psi]
\arrow|r|[b`c;\psi \tens \psi]
\arrow|r|[c`d;1 \tens a \tens 1]
\arrow|r|[d`e;\phi \tens \phi]
\arrow|r|[e`f;G\phi \tens G\phi]
\arrow|l|[a`g;\phi]
\arrow|l|[g`h;G\phi]
\arrow|l|[h`i;GF(1 \tens y \tens 1)]
\arrow|l|[i`j;G\psi]
\arrow|l|[j`f;\psi]
\arrow|m|/{@{>}@/_2em/}/[b`k;\phi]
\arrow|m|/{@{>}@/^2em/}/[g`k;G(\psi \tens \psi)]
\arrow|m|[k`l;G(1 \tens z \tens 1)]
\arrow|m|/{@{>}@/_2em/}/[l`e;\psi]
\arrow|m|/{@{>}@/^2em/}/[l`j;G(\phi \tens \phi)]

\efig \]

\caption{Proof of Proposition \ref{pfc}}
\label{lozenge}
\end{figure}

\noindent The diamonds commute by naturality of $\phi$ and $\psi$ and the left and right pentagons commute by 
prebimonoidality of $F$ and $G$, respectively. \end{proof}

\begin{proposition} \label{prop310}
If $F$ is separable Frobenius then it is prebimonoidal relative to $y$ and $z=y^{F}$.
\end{proposition}
\begin{proof} The proof is contained in Figure \nolinebreak \ref{crystal}.
\begin{figure}[htbp]
\[ \bfig

\hscalefactor{1.2}
\vscalefactor{0.9}
\node a(0,+2000)[F(Tu \tens Tv) \tens F(Tw \tens Tx)]
\node b(+500,+1500)[F(Tu \tens Tv \tens Tw \tens Tx)]
\node c(+500,+500)[F(Tu \tens Tv \tens Tw \tens Tx)]
\node d(+500,-750)[F(Tu \tens Tw \tens Tv \tens Tx)]
\node e(+500,-1750)[F(Tu \tens Tw \tens Tv \tens Tx)]
\node f(0,-2250)[F(Tu \tens Tw) \tens F(Tv \tens Tx)]
\node g(-500,+1500)[FTu \tens FTv \tens F(Tw \tens Tx)]
\node h(-1000,+1000)[FTu \tens FTv \tens FTw \tens FTx]
\node i(-500,+500)[FTu \tens F(Tv \tens Tw) \tens FTx]
\node j(-500,-750)[FTu \tens F(Tw \tens Tv) \tens FTx]
\node k(-1000,-1250)[FTu \tens FTw \tens FTv \tens FTx]
\node l(-500,-1750)[F(Tu \tens Tw) \tens FTv \tens FTx]
\node m(0,+1000)[FTu \tens F(Tv \tens Tw \tens Tx)]
\node n(0,-1250)[F(Tu \tens Tw \tens Tv) \tens FTx]
\node o(0,0)[F(Tu \tens Tv \tens Tw) \tens FTx]

\arrow|r|[a`b;\phi]
\arrow/=/[b`c;]
\arrow|r|[c`d;F(1 \tens y \tens 1)]
\arrow/=/[d`e;]
\arrow|r|[e`f;\psi]
\arrow|m|[a`g;\psi \tens 1]
\arrow|m|[g`h;1 \tens 1 \tens \psi]
\arrow|m|[h`i;1 \tens \phi \tens 1]
\arrow|m|[i`j;1 \tens Fy \tens 1]
\arrow|m|[j`k;1 \tens \psi \tens 1]
\arrow|m|[k`l;\phi \tens 1 \tens 1]
\arrow|m|[l`f;1 \tens \phi]
\arrow|m|[g`m;1 \tens \phi]
\arrow|m|[b`m;\psi]
\arrow|m|[m`c;\phi]
\arrow|m|[m`i;1 \tens \psi]
\arrow|m|[j`n;\phi \tens 1]
\arrow|m|[d`n;\psi]
\arrow|m|[n`e;\phi]
\arrow|m|[n`l;\psi \tens 1]
\arrow|l|[h`k;1 \tens y^F \tens 1]
\arrow|m|[i`o;\phi \tens 1]
\arrow|m|[c`o;\psi]
\arrow|m|[o`n;F(1 \tens y) \tens 1]
\arrow|l|/{@{>}@/_6em/}/[a`h;\psi \tens \psi]
\arrow|l|/{@{>}@/_6em/}/[k`f;\phi \tens \phi]

\efig \]
\caption{Proof of Proposition \ref{prop310}}
\label{crystal}
\end{figure}

\noindent The five diamonds commute since $F$ is Frobenius, and the two right-hand triangles commute since $F$ is separable. 
The rhombus commutes by definition
of $y^F$, the parallelograms by naturality of $\phi$ and $\psi$, and the two irregular cells are trivial.
\end{proof}

\begin{proposition} \label{prop311}
A strong monoidal functor between braided monoidal categories is prebimonoidal if and only if it is braided.
\end{proposition}
\begin{proof}
As noted above, strong monoidal functors are separable Frobenius, and strong monoidal functors are braided precisely when $c_{A,B}^F = c_{FA,FB}$,
so Proposition \nolinebreak \ref{prop310} establishes ``if''. Conversely, suppose that $F$ is prebimonoidal with respect to the two braidings, and consider
the commutative diagram in Figure \nolinebreak \ref{hexes}.

\begin{figure}[htbp]
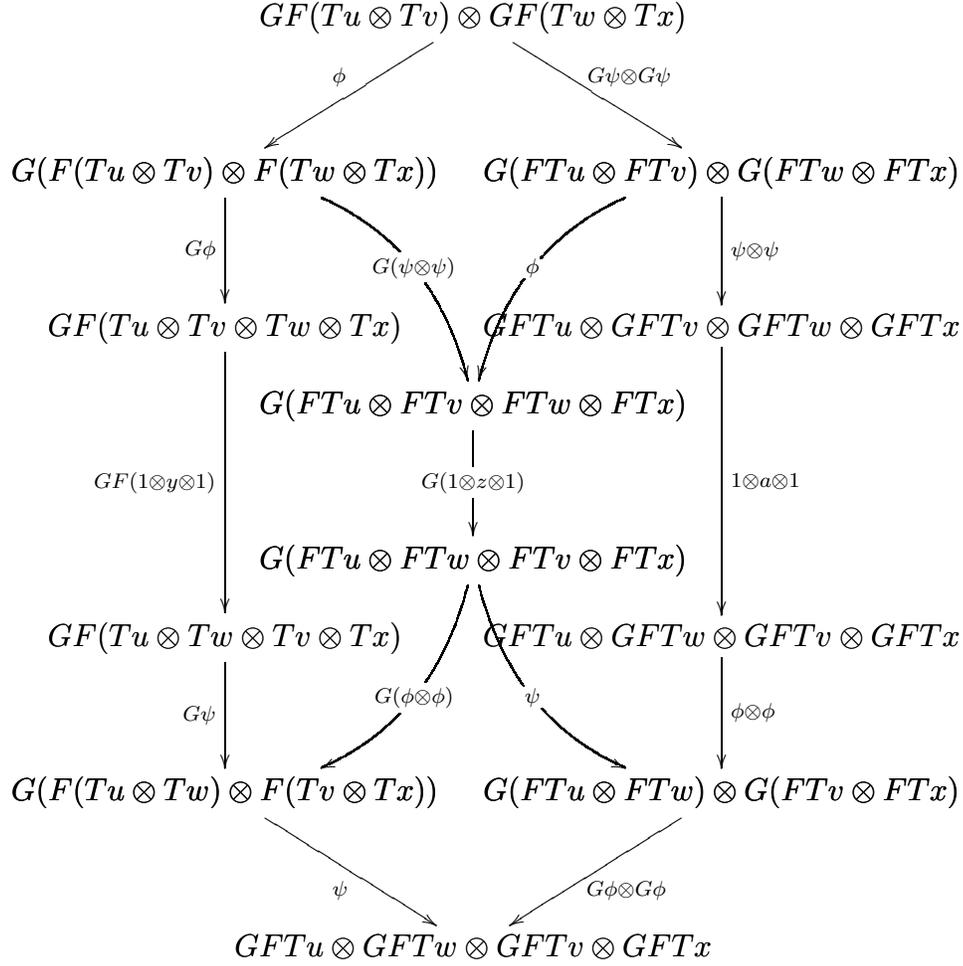

\[ \bfig

\hscalefactor{1.2}
\node aa(0,+1000)[F(I \tens A) \tens F(B \tens I)]
\node bb(+500,+500)[F(I \tens A \tens B \tens I)]
\node cc(+500,0)[F(I \tens B \tens A \tens I)]
\node dd(0,-500)[F(I \tens B) \tens F(A \tens I)]
\node ee(-500,+500)[FI \tens FA \tens FB \tens FI]
\node ff(-500,0)[FI \tens FB \tens FA \tens FI]
\node a(0,+1500)[FA \tens FB]
\node b(+1000,+1000)[F(A \tens B)]
\node c(+1000,-500)[F(B \tens A)]
\node d(0,-1000)[FB \tens FA]
\node e(-1000,+1000)[I \tens FA \tens FB \tens I]
\node f(-1000,-500)[I \tens FB \tens FA \tens I]


\arrow|r|[a`b;\phi]
\arrow|r|[b`c;Fc]
\arrow|r|[c`d;\psi]
\arrow/=/[a`e;]
\arrow|l|[e`f;1 \tens c \tens 1]
\arrow/=/[f`d;]
\arrow|m|[aa`bb;\phi]
\arrow|m|[bb`cc;F(1 \tens c \tens 1)]
\arrow|m|[cc`dd;\psi]
\arrow|m|[aa`ee;\psi \tens \psi]
\arrow|m|[ee`ff;1 \tens c \tens 1]
\arrow|m|[ff`dd;\phi \tens \phi]
\arrow/=/[a`aa;]
\arrow/=/[b`bb;]
\arrow/=/[c`cc;]
\arrow/=/[d`dd;]
\arrow|m|[ee`e;\psi_0 \tens 1 \tens 1 \tens \psi_0]
\arrow|m|[f`ff;\phi_0 \tens 1 \tens 1 \tens \phi_0]

\efig \]
\caption{Proof of Proposition \ref{prop311}}
\label{hexes}
\end{figure}

The middle cell commutes since $F$ is prebimonoidal, the bottom left since $\phi$ is monoidal, and the top left since $\psi$ is opmonoidal.
The three right-hand cells commute by definition, and, noting that $\phi_0$ and $\psi_0$ are both natural and mutually inverse, the 
left-hand cell does so also.
Hence, the full diagram shows that $c_{A,B}^F = c_{FA,FB}$, as desired. \end{proof}


\bibliographystyle{plain}

\vspace{3cm}

\noindent M. B. McCurdy, Macquarie University, NSW 2109, Australia

\noindent R. Street, Macquarie University, NSW 2109, Australia

\noindent The authors are grateful for the support of the Australian Research Council and Macquarie University.

\end{document}